\documentclass[leqno]{amsart}

\usepackage[all]{xy}
\usepackage{amsmath, amscd,amsthm}

\swapnumbers

\theoremstyle{plain}
\newtheorem{theorem}{Theorem}[section]
\newtheorem{proposition}[theorem]{Proposition}
\newtheorem{lemma}[theorem]{Lemma}

\newtheorem{cory}[theorem]{Corollary}
\newtheorem{prop}[theorem]{Proposition}

\newtheorem{remarks}[theorem]{Remarks}

\theoremstyle{definition}

\newtheorem{definition}[theorem]{Definition}

\newtheorem{remark}[theorem]{Remark}

\numberwithin{equation}{section}

\newcommand{\secref}[1]{\S\ref{#1}}
\newcommand{\thmref}[1]{Theorem~\ref{#1}}
\newcommand{\propref}[1]{Proposition~\ref{#1}}
\newcommand{\lemref}[1]{Lemma~\ref{#1}}
\newcommand{\corref}[1]{Corollary~\ref{#1}}

\newcommand{\remref}[1]{Remark~\ref{#1}}

\newcommand{\defref}[1]{Definition~\ref{#1}}

\def\m{\medskip}

\def\ov{\overline}

\def\cat{{\rm cat}}
\def\cl{{\rm cl}}
\def\dim{{\rm dim}}
\def\wgt{\rm wgt}
\def\cwgt{\rm cwgt}
\def\swgt{\rm swgt}

\def\PP{\mathbb P}
\def\RR{\mathbb R}
\def\ZZ{\mathbb Z}
\def\rp2{\RR\PP^2}
\def\Hom{{\rm Hom}}
\def\connsum{\,\#\,}

\begin{document}

\title{Detecting Elements and Lusternik--Schnirelmann
Category of 3-Manifolds}

\author{John Oprea}
\address{Department of Mathematics\\
         Cleveland State University\\
         Cleveland Ohio 44115  U.S.A.}
   
\email{oprea@math.csuohio.edu}

\author{Yuli Rudyak}
\address{Department of Mathematics\\
         University of Florida\\
         358 Little Hall\\
         PO Box 118105\\
         Gainesville, FL 32611-8105  U.S.A.}
   
\email{rudyak@math.ufl.edu}

\date{\today}

\begin{abstract}
In this paper, we give a new simplified calculation
of the Lusternik-Schnirelmann category of closed
$3$-manifolds. We also describe when $3$-manifolds
have detecting elements and prove that $3$-manifolds
satisfy the equality of the Ganea conjecture. 

\end{abstract}


\maketitle

%\pagestyle{myheadings}
%\markboth{John Oprea and Yuli Rudyak}
%{Detecting Elements and LS-Category of $3$-Manifolds}

\section{Introduction}
The \emph{Lusternik--Schnirelmann category} of a space $X$,
denoted $\cat(X)$, is defined to be the minimal integer $k$
such that there exists an open covering $\{A_0, \ldots, A_k\}$
of $X$ with each $A_i$ contractible to a point in $X$.
Category, while easy to define, is notoriously difficult to
compute in general. In particular, except for $K(\pi,1)$'s,
it cannot be expected that the category of a space is
determined by its fundamental group. In \cite{GoGo},
however, the following interesting result was proved.

\begin{theorem}\label{thm:3mancat}
Let $M^3$ be a closed $3$-dimensional manifold. Then
$$\cat(M) =
\begin{cases}
1 & {\rm if}\quad \pi_1(M) = \{1\} \\
2 & {\rm if}\quad \pi_1(M) \ {\rm is\ free} \\
3 & {\rm otherwise}
\end{cases}.
$$
\end{theorem}

In this paper, we will give a somewhat simplified
proof of this theorem using the relatively new
approximating invariant for category, \emph{category
weight}. Throughout, we use only basic results about
$3$-manifolds found, for instance, in \cite{H}.
\emph{But we shall also do more}. We will prove that
most 3-manifolds possess a \emph{detecting element};
that is, an element whose category weight is equal to
the category of $M$ (see \cite{R3}). It is known that
a \emph{detectable space} (i.e., a space possessing
detecting elements) has some special properties which allow
solutions of certain well-known problems (\cite{R3}).
For example, from the existence of detecting elements,
we prove that closed 3-manifolds satisfy the
Ganea conjecture.

\begin{cory}
For every closed $3$-manifold $M$,
$$\cat(M\times S^n) = \cat(M)+1 .$$

\end{cory}

\emph{This result is not obtainable from knowing
the category alone}, so the detecting element
approach is a significant embellishment of
\thmref{thm:3mancat}. Another well-known problem
is the relationship between degree $1$ maps of
manifolds and LS-category. For closed, $3$-manifolds,
we have 

\begin{cory}
Let $f\colon M \to N$ be a degree $1$ map of oriented
$3$-manifolds. Then $\cat M \ge \cat f = \cat N$.
\end{cory}

We now turn to the fundamentals of $3$-manifolds.

\section{Preliminaries on $3$-Manifolds}

\begin{definition}\label{def:irred-prime}
A $3$-manifold $M$ is {\emph{irreducible}} if every
embedded two-sphere $S^2 \hookrightarrow M$ bounds an
embedded disk $D^3 \hookrightarrow M$.

A $3$-manifold $M$ is \emph{prime} if $M = P \connsum Q$ implies
that either $P = S^3$ or $Q = S^3$. Here, ``$=$'' denotes
diffeomorphism and $\connsum$ is the connected sum.
\end{definition}

The following two results clarify the relation between
prime and irreducible manifolds.
 
\begin{lemma}\label{lem:irred-prime}
If $M^3$ is irreducible, then it is prime.
\end{lemma}

\begin{proof}
Suppose $M$ is irreducible. In order to split $M$ as
$M = P \connsum Q$, there must be an embedded $S^2$ which separates
$M$ into two components (i.e. $P - D^3$ and $Q - D^3$). But
any such $S^2$ bounds an embedded disk $D^3$ by irreducibility,
so $M$ can only split as $M = M' \connsum S^3$ (since $S^3 - D^3$
is a disk $D^3$). This says that $M$ is prime.
\end{proof}

\begin{lemma}\label{lem:bundle}
If $M$ is a prime 3-manifold and $M$ is not irreducible,
then $M$ is the total space of a 2-sphere bundle over $S^1$.
\end{lemma}

\begin{proof}
See \cite[Lemma 3.13]{H}
\end{proof}

\noindent The fundamental structural result about $3$-manifolds
is the following 

\begin{theorem}[{Prime Decomposition}]
\label{thm:prim-dec} 
 A $3$-manifold $M$ may
be written as 
$$
M = M_1 \connsum M_2 \connsum \ldots \connsum M_k,
$$
where each $M_j$ is prime. Furthermore, such a prime
decomposition is unique up to re-arrangement of summands.
\end{theorem}

\begin{proof}
See \cite[Theorems 3.15 and 3.21]{H}
\end{proof}

The \emph{Sphere theorem} says that, for an orientable
$3$-manifold $M$, $\pi_2(M) \not = 0$ implies that
some element of $\pi_2(M)$ is represented by an
embedding $S^2 \hookrightarrow M$. We will require
the following generalization.

\begin{theorem}[The Projective Plane Theorem]
\label{proj-plane}
Let $M$ be a 3-manifold with \linebreak
$\pi_2(M)\ne 0$. Then
there exists a map $g\colon S^2 \to M$ with the following
properties.
\begin{enumerate}
\item The map $g$ is not null-homotopic.
\item The map $g\colon S^2 \to g(S^2)$ is a covering map.
\item $g(S^2)$ is a $2$-sided submanifold ($2$-sphere or
projective plane) in $M$.
\end{enumerate}
\end{theorem}

\begin{proof}
See \cite[Theorem 4.12]{H}.
\end{proof}

\noindent With these preliminaries, we can prove the
folowing important characterization. 

\begin{prop}\label{irredchar} 
Let $M$ be a closed $3$-manifold. Then, 
\begin{enumerate}
\item  If $\pi=\pi_1(M)$ is infinite and $\pi_2(M)=0$,
then $M = K(\pi,1)$.
\item If $\pi_1(M)$ is finite, then the universal covering
of $M$ is a homotopy $3$-sphere and $M$ is orientable.
\end{enumerate}
\end{prop}

\begin{proof} For (1), assume that $\pi_1(M)$ is infinite.
Let $p\colon \widetilde M \to M$ be the universal covering
of $M$. Since $\pi_2(M)=0$, we conclude that 
$H_2(\widetilde M)=0$. Since $\pi_1(M)$ is infinite, we
conclude that $\widetilde M$ is not compact, and therefore
$H_3(\widetilde M)=0$. Hence, $\widetilde M$ is acyclic.
Moreover, $\widetilde M$ is simply-connected, and, by the
Whitehead theorem, it is therefore contractible. Hence,
$M=K(\pi_1(M),1)$. 

\m For (2), assume that $\pi_1(M)$ is finite. Then the
universal cover $\widetilde M$ of $M$ is a closed simply
connected manifold. So, by Poincar\'e duality,
$H_2(\widetilde M)=0$, and hence, by the Hurewicz theorem, 
$\pi_2(\widetilde M)=0$. Thus, $\pi_2(M)=0$. Furthermore, 
$$
H_3(\widetilde M) = \ZZ = \pi_3(\widetilde M),
$$
again by the Hurewicz theorem. Therefore, the generator of
$\pi_3(\widetilde M)=\ZZ$ provides a degree $1$ map $S^3
\to \widetilde M$ (i.e. an isomorphism on $H_3$). Since
$\widetilde M$ and $S^3$ are simply connected, the Whitehead
theorem implies that $\widetilde M \simeq S^3$.

To see that $M$ is orientable, we simply note that each
$g \in \pi_1(M)$, thought of as a covering transformation
on the orientable manifold $\widetilde M$, acts to preserve
orientation. This is seen by supposing the opposite;
namely, that $g$ reverses orientation. Now, because
$\widetilde M \simeq S^3$, homotopy classes of maps
$\widetilde M \to \widetilde M$ are classified by degree.
Since $g$ is a homeomorphism which reverses orientation,
its degree is $-1$. But then the Lefschetz number of $g$ is
$L(g) = 2$, implying the existence of a fixed point and
contradicting the fact that $g$ is a covering transformation.
Hence, all covering transformations preserve orientation,
so $M = \widetilde M/\pi_1(M)$ is orientable.
\end{proof}

\noindent These are the only ingredients from $3$-manifold
theory that we shall need. In the next section, we introduce
the main technical tool, the approximating invariant category
weight.

\section{Category Weight and Detecting Elements}

\begin{definition}[\cite{BG,Fe,F}]
\label{def:cat-f}
Let $f\colon X \to Y$ be a map of finite $CW$-spaces. The 
\emph{Lusternik--Schnirelmann  category of $f$}, denoted
$\cat(f)$, is defined to be the minimal integer $k$ such
that there exists an open covering $\{A_0, \ldots, A_k\}$
of $X$ with the property that each of the restrictions
$f|A_i\colon A_i \to Y$, $i=0,1, \ldots, k$ is null-homotopic.
\end{definition}

\noindent Clearly, $\cat(X) = \cat(1_X)$ and. Also, it is
easy to see that $\cat(f) \leq \cat(X)$ since $f$ is
null-homotopic on any subset which is contractible in $X$.

\begin{definition}\label{def:swgt} 
The \emph{category weight} of a non-zero cohomology class
$u \in H^*(X; R)$ (for some, possibly local, coefficient
ring $R$) is defined by 
$$
\wgt(u) \geq k {\rm\ if\ and\ only\ if\ } \phi^*(u)=0 {\rm\
for\ any\ } \phi\colon A \to X {\rm\ with\ } \cat(\phi) < k .$$
\end{definition}

\begin{remarks}\label{rem:credits}\rm
1. The idea of category weight was suggested by Fadell and
Husseini (see \cite{FH}). In fact, they considered an invariant
similar to our $\wgt$ (denoted in \cite{FH} by $\cwgt$), 
but where the defining maps $\phi\colon A \to X$ were
required to be inclusions instead of general maps. Because
of this, $\cwgt$ was not a homotopy invariant, and this made
it a delicate quantity in homotopy calculations. Rudyak in
\cite{R2, R3} and Strom in \cite{S} suggested the homotopy
invariant version of category weight as defined in
\defref{def:swgt}. Rudyak called it \emph{strict} category
weight (using the notation $\swgt (u)$) and Strom called it
\emph{essential} category weight (using the notation $E(u)$).
At the Mt. Holyoke conference for which these proceedings
are a record, both creators agreed to adopt the notation
$\wgt$ and call it simply \emph{category weight}.

\m 
2. In fact, one can define category weight for $u\in F^*(X)$
where $F$ is a suitable functor on the category of topological
spaces (e.g. $F(X)=[X,Y]$ or $F$ is an arbitrary cohomology
theory), see \cite{R2, R3, S}. However, \defref{def:swgt} is
enough for our goals here.

\m
3. There is an alternative definition of category weight
which is actually more useful than the one given in
\defref{def:swgt}. Recall that the Ganea fibration
$p_j \colon G_j(X) \to X$ is defined inductively starting with
the path fibration $p_0 \colon PX = G_0(X) \to X$ having
fibre $\Omega X$. Then given the fibration $p_i\colon G_i(X)
\to X$ with fibre $F_i = *^{(i+1)}\Omega X$, the fibration
$p_{i+1}$ is constructed by taking the cofibre $Z$ of the
inclusion $F_i \to G_i(X)$ and extending $p_i$ to a map
$Z \to X$ (which is possible since the composition
$F_i \to G_i(X) \stackrel{p_i}{\to} X$ is null-homotopic.
Finally, convert the map $Z \to X$ to a fibration $p_{i+1}
\colon G_{i+1}(X) \simeq Z \to X$. Then it is known that
$\cat(X) = k$ if and only if $k$ is the least integer
such that $p_k\colon G_k(X) \to X$ has a section, \cite{G, Sv}.
It can also be shown that, for a cohomology class $u \in H^*(X;R)$,
$\wgt(u) = k$ if and only if $k$ is the greatest integer
such that $p_{k-1}^*(u) = 0$, \cite{R3, S}. We shall use this below in
giving a proof of \propref{prop:swgtprops} (4).
\end{remarks}

\begin{proposition}[\cite{R3,S}]
\label{prop:swgtprops}
Category weight has the following properties.
\begin{enumerate}
\item $1\le \wgt(u) \leq \cat(X)$, for all $u \in \widetilde
H^*(X;R), u\ne 0$.
\vskip3pt
\item For every $f\colon Y \to X$ and $u\in H^*(X;R)$ with
$f^*(u)\not = 0$ we have
$\cat(f) \geq \wgt(u)$ and $\wgt(f^*(u)) \geq \wgt(u)$.
\vskip3pt
\item $\wgt(u\cup v) \geq \wgt(u) + \wgt(v)$. 
\vskip3pt
\item For every $u \in H^s(K(\pi,1);R)$, $u\ne 0$, we have
$\wgt(u)\geq s$.
\vskip3pt
\end{enumerate}
\end{proposition}

\begin{proof} We will only prove (4) since the other results
are proven in the references cited. If $X = K(\pi,1)$, then
$\Omega X$ has the homotopy type of a discrete set of points
and, consequently, $F_1 = \Omega X * \Omega X$ is, up to
homotopy, a wedge of circles. Also, $G_0(X) = PX
\simeq *$, so the cofibre of $\Omega X \to G_0(X)$ has the
type of a wedge of circles. Then $G_1(X)$ has the homotopy
type of a $1$-dimensional space. Similarly, it is easy to
see that $G_k(X)$ has the homotopy type of a $k$-dimensional
space. If $u \in H^s(K(\pi,1);R)$, then $p_{s-1}^*(u) = 0$
since $G_{s-1}(X)$ is $s$-dimensional. By the equivalent
definition of $\wgt$ given in \remref{rem:credits} (3), we
see that $\wgt(u) \geq s$.
\end{proof}

\begin{definition}
We say that $u\in H^*(X;R)$ is a \emph{detecting element}
for $X$ if $\wgt(u) =\cat(X)$. We say that a space $X$ is
\emph{detectable} if it possesses a detecting element.
\end{definition}

Recall that the \emph{cup-length} of a space $X$ with respect
to a ring $R$ is defined as
$$
\cl_R(X)=\max\{k\bigm| u_1\cup \cdots \cup u_k \ne 0\
\text{ for some } u\in \widetilde H^*(X;R)\}.
$$

\begin{lemma}\label{cup-length}
If $\cat(X) = \cl _R(X)$ for some ring $R$ then the space $X$
is detectable. 
\end{lemma}

\begin{proof}
It is well known that $\cat(X) \ge \cl_R(X)$ for every $R$.
Now, let $\cat(X) =k$ and suppose that there are $u_1, \ldots,
u_k\in \widetilde H^*(X;R)$ with $u_1\cup \cdots \cup u_k \ne 0$.
Then, using the first and third properties of
\propref{prop:swgtprops}, we conclude that $\wgt(u_1\cup \cdots
\cup u_k)=k$. Thus, $u_1\cup \cdots \cup u_k$ is a detecting
element for $X$.
\end{proof}

\section{Basic Special Cases}\label{sec:basic}

First, recall that $\cat(X)  \le \dim(X)$ for every connected
$CW$-space $X$. In particular, $\cat(M) \le 3$ for every (connected)
$3$-manifold $M$. We also notice that, by \lemref{cup-length}, a
space $X$ is detectable whenever $\cat(X) = \cl_{R}(X)$ for some $R$.
Here is a first step in understanding the category of $3$-manifolds.

\begin{prop}\label{finite}
If $M$ is a 3-manifold with finite fundamental group of order $d>1$,
then $\cat(M)=3$, and every non-zero element of $H^3(M;\ZZ/d)$ is
a detecting element for $M$. Moreover, if $d$ is even, then every
non-zero element of $H^3(M;\ZZ/2)$ is a detecting element for $M$
as well.
\end{prop}

\begin{proof} Since $\pi_1(M)$ is finite, $\pi_2(M)=0$ because, by 
\propref{irredchar}, the universal cover is a homotopy sphere.
Hence. there is the Hopf exact sequence
$$
\CD
\pi_3(M) @>h>> H_3(M) @>q>> H_3(\pi) \to 0
\endCD
$$
where $h$ is the Hurewicz homomorphism (e.g. see
\cite[Theorem II.5.2]{Br}. Since, by \propref{irredchar},
the $d$-fold universal covering $\widetilde M \to M$ is a
$d$-sheeted covering, $M$ is orientable and $\widetilde M$
is a homotopy sphere, we conclude that $h$ has the form
$$
\pi_3(M)=\ZZ \to \ZZ=H_3(M),\quad a \mapsto d\cdot a.
$$
Hence, $H_3(\pi)=\ZZ/d$. Also consider the induced
homomorphism ${\rm Hom}(H_3(\pi);\ZZ/d) \to
{\rm Hom}(H_3(M); \ZZ/d)$. It is certainly injective
since $H_3(M) \to H_3(\pi)$ is surjective. However,
it is also true that, for any $\phi \in
{\rm Hom}(H_3(M); \ZZ/d)$, ${\rm Im}(h) = d\ZZ \subseteq
{\rm Ker}(\phi)$, so there exists $\bar \phi \in
{\rm Hom}(H_3(\pi);\ZZ/d)$ with $\bar\phi \mapsto \phi$.
Thus, we have an isomorphism ${\rm Hom}(H_3(\pi);\ZZ/d)
\stackrel{\cong}{\to} {\rm Hom}(H_3(M); \ZZ/d)$.

Now consider the diagram
$$
\CD
H^3(\pi; \ZZ/d) @>q^*>> H^3(M;\ZZ/d)\\
@VVV @VVV\\
\Hom(H_3(\pi);\ZZ/d) @>q^*>> \Hom(H_3(M); \ZZ/d).
\endCD
$$
By \propref{prop:swgtprops}, (4), a non-zero element of
$H^3(\pi;\ZZ/d)$ has category weight at least $3$. The right
arrow is an isomorphism because $H_2(M)$ is free
abelian since $M$ is orientable. The
bottom arrow is an isomorphism by the argument above. Finally,
the left arrow is a surjection by the Universal Coefficient
Formula. Therefore, the top arrow is a surjection as well.
In particular, by \propref{prop:swgtprops} (2), every
non-zero element of $H^3(M;\ZZ/d)$ has category weight at
least $3$. But $\cat(M) \le \dim(M) =3$, so $\cat(M)=3$,
and every non-zero element of $H^3(M;\ZZ/d)$ is a detecting
element for $M$.
\end{proof}

\begin{remark} Using the approach as in \propref{finite},
it is also possible to prove the following result originally
due to Krasnoselski \cite{Kra} and, in fact, re-proved
in \cite{GoGo}:

For a free action of the finite group $G$ on a homotopy
sphere $S \simeq S^{2n+1}$,
$$
\cat(S/G) = 2n+1 = \dim(S/G).
$$
\end{remark}

Here is another basic result which follows from the
characterization of prime non-irreducible $3$-manifolds.

\begin{prop}\label{cat=2}
Let $M$ be a prime $3$-manifold which is not irreducible.
Then $\cat(M) = 2 = \cl_{\ZZ/2}(M)$, and $M$ is detectable. 
\end{prop}

\begin{proof} In view of \lemref{lem:bundle}, $M$ is the
total space of a $2$-sphere bundle over $S^1$. So, $M$ is
either $S^1 \times S^2$ or the mapping torus of the map 
$$
r\colon S^2 \to S^2, \quad r(x)=-x
$$
where $S^2$ is regarded as the set of unit vectors in
$\RR^3$. It is easy to see that, in both of the cases,
$M = (S^1 \vee S^2)\cup e^3$ where $e^3$ is a $3$-cell
attached to the wedge $S^1 \vee S^2$. Thus, because a
wedge of spheres has category one and a mapping cone
can increase category by at most one, we obtain
$\cat(M) \le 2$.  

Futhermore, because $\pi_1(M)=\ZZ$, we conclude that
$H_1(M;\ZZ/2)=\ZZ/2$. So, because of Poincar\'e duality
(with $\ZZ/2$-coefficients), we have $\cl_{\ZZ/2}(M)\ge 2$.
Thus, $\cl_{\ZZ/2}(M)=2=\cat(M)$, and $M$ is detectable.
\end{proof}

The next two results treat the case of infinite fundamental
group, excluding the $S^2$-bundles over $S^1$.
 
\begin{prop}\label{infinite}
If $M$ is a $3$-manifold with $\pi_1(M)$ infinite and
$\pi_2(M)=0$, then $\cat(M)=3$ and $M$ is detectable.  
\end{prop}

\begin{proof}  By \propref{irredchar}, $M=K(\pi_1(M),1)$,
so, by \propref{prop:swgtprops}, every non-zero element of
$H^3(M;R)$ has category weight $3$. (Notice that, for
example, $H^3(M;\ZZ/2)\ne 0$). Thus, because
$\cat(M) \le \dim(M) = 3$, each of these elements is
a detecting element.
\end{proof}

\begin{prop}\label{non-or}
If $M$ is an irreducible $3$-manifold with $\pi_1(M)$
infinite and \linebreak
$\pi_2(M) \ne 0$, then $\cat(M) = 3 =
\cl_{\ZZ/2}(M)$. In particular, $M$ is detectable. 
Furthermore, $M$ is non-orientable. 
\end{prop}

\begin{proof}
Consider a map $g\colon S^2 \to M$ as in \thmref{proj-plane}.
Since $M$ is irreducible, we conclude that $g(S^2)$ is a
$2$-sided projective plane in $M$. Let $i\colon \rp2\to M$
be the corresponding embedding, and let $[\rp2]\in 
H_2(\rp2;\ZZ/2)$ denote the fundamental class modulo $2$
of $\rp2$.

Let $w_k$ and $\ov w_k$ denote the $k$-th Stiefel--Whitney
class of $M$ and $\rp2$, respectively. Since the $1$-dimensional
normal bundle of $i$ is trivial, we conclude that $i^*w_k
= \ov w_k$. We can now compute the Kronecker products
$$
\langle w_2, i_*[\rp2] \rangle = \langle i^*w_2, [\rp2] \rangle
= \langle \ov w_2, [\rp2] \rangle = 1,
$$  
and so $i_*[\rp2] \ne 0 \in H_2(M;\ZZ/2)$. Now, since
$\langle \ov w_1^2, [\rp2] \rangle = 1$, we conclude that
$i^*w_1^2 = \ov w_1^2 \ne 0$, and so $w_1^2\ne 0$. 
So, by Poincar\'e duality, there exists $x\in H^1(M;\ZZ/2)$
with $x w_1^2 \ne 0$. Thus, $\cl_{\ZZ/2}(M)=3$.
\end{proof}

We also need the following fact which, in a sense,
is a converse of \lemref{lem:bundle}.

\begin{cory}\label{fundgroup-free}
If $M$ is a closed $3$-manifold with non-trivial free
fundamental group, then $M$ is not irreducible.
\end{cory}

\begin{proof} 
Notice that $\pi_2(M)\ne 0$. Indeed, if $\pi_2(M)=0$
then, by \propref{irredchar} and the hypothesis that
$\pi_1(M)$ is free, 
$$
M = K(\pi_1(M),1) = \vee S^1.
$$
But this is wrong since a wedge of circles has vanishing
homology above degree $1$ for any coefficients.

Now, if $M$ is irreducible then, by \propref{non-or},
$\cl_{\ZZ/2}(M)=3$. But this is impossible. Indeed, let 
$f\colon M \to K(\pi_1(M),1) = \vee S^1$ be a map which
induces an isomorphism of fundamental groups. Then 
$$
f^*\colon H^1(K(\pi_1(M), 1);\ZZ/2) \to H^1(M;\ZZ/2)
$$
is an isomorphism. Thus, $x \cup y=0$ for all
$x,y \in H^1(M;\ZZ/2)$, and so $\cl_{\ZZ/2}(M) < 3$.
This is a contradiction.   
\end{proof}

\begin{remark}
If $\pi_1(M)=\ZZ$ then $M = P \connsum\Sigma$ where $\Sigma$
is a homotopy sphere and $P$ is prime. So, $\pi_1(P)=\ZZ$.
But $P$ is not irreducible by \corref{fundgroup-free},
so, because of \lemref{lem:bundle}, $\pi_2(P)=\ZZ$.
In other words, $\pi_2(M)=\ZZ$ whenever $\pi_1(M)=\ZZ$. 
Actually, the following general fact holds: for every
closed 3-manifold $M$, the group $\pi_1(M)$ completely
determines $\pi_2(M)$, see e.g. \cite{R1}.
\end{remark}

\section{Detectability of $3$-Manifolds}

\begin{proposition}\label{prop:cat2free}
Let $M^3$ be a closed $3$-manifold with $\pi_1(M)$
free and non-trivial. Then $\cat(M) = 2$, and $M$ is detectable.
\end{proposition}

\begin{proof}
Write $M = M_1 \connsum \ldots \connsum M_k$ with each $M_j$
prime. Because $\pi_1(M) = \pi_1(M_1) * \ldots *
\pi_1(M_k)$ is free, each $\pi^j = \pi_1(M_j)$ must
be free (where we agree that the trivial group is free).
If $M_j$ is irreducible with $\pi^j \not = \{1\}$, then
this contradicts \corref{fundgroup-free}. Therefore,
all such $M_j$ are non-irreducible primes; that is,
the $M_j$ are the manifolds considered in \propref{cat=2}.
Because of \lemref{lem:bundle}, these are the total
spaces of $S^2$-bundles over $S^1$. There are only two
such manifolds: one orientable and one non-orientable,
and we denote both of them by $S^1\propto S^2$.
Of course, the $M_j$ with $\pi^j = \{1\}$ are homotopy
spheres $\Sigma_j$. The key point now is that, for
$M = P \connsum Q$ with $P = \connsum_k (S^1 \propto S^2)$ and
$Q = \connsum_j \Sigma_j$, $M - D^3$ deformation retracts
onto the $2$-skeleton $\vee_k (S^1 \vee S^2)$. Because
of \propref{cat=2}, $\cat(S^1 \propto S^2) = 2$. This
handles the ``trivial'' case where the connected sum
degenerates to a single summand. Now suppose
$M = \connsum_j M_j = P \connsum Q$, where $M_j$ is either a
homotopy sphere or $S^1 \propto S^2$ and $P =
\connsum_{j_t} M_{j_t}$, $Q = \connsum_{j_s} M_{j_s}$ arbitrarily
split $M$. If we remove a disk from a $3$-manifold
$N$, then the inclusion $S^2 \hookrightarrow N - D^3$
is the inclusion of a subcomplex; so therefore a
cofibration. Thus, the pushout diagram
$$\xymatrix{
S^2 \ar[r] \ar[d] & P - D^3 \ar[d] \\
Q - D^3 \ar[r] & P \connsum Q = M }
$$

\noindent is a homotopy pushout as well. But then we
may apply the standard estimate for the category
of a double mapping cylinder (see \cite{Har}) to obtain
\begin{align*}
\cat(M) & \leq \cat(S^2) + \max\{\cat(P - D^3),
\cat(Q - D^3)\} \\
& = 1 + \max\{\cat(\vee_{j_t}(S^1 \vee S^2)),
\cat(\vee_{j_s}(S^1 \vee S^2))\} \\
& = 1 + 1 \\
& = 2.
\end{align*}

\noindent Of course, cup-length then shows that $\cat(M)
= 2$ and this completes the proof.
\end{proof}

\begin{theorem}\label{thm:non-free}
Let $M$ be a $3$-manifold whose fundamental group
is non-trivial and not a free group. Then
$\cat(M) = 3$. Further, $M$ is detectable unless
it is non-orientable of the form $P \connsum Q$, where
$P$ is non-orientable and $Q$ is prime with odd
torsion. Also, in the last case, the orientable
double cover of $M$ has category $3$.
\end{theorem}

\begin{proof}
The case of finite $\pi_1$ is considered in
\propref{finite}. So, we assume that $\pi_1(M)$ is
infinite. We represent $M$ as a connected sum
$M = N\connsum P$, where $P$ is prime and $\pi_1(P) \ne
\{1\}$. Furthermore, we can always assume that 
$\pi_1(P)\ne \ZZ$, and therefore $P$ is irreducible
in view of \corref{fundgroup-free}. Now, because of
the results of \secref{sec:basic}, $P$ posseses a
detecting element $u\in H^3(P;R)$ for suitable $R$.

Now suppose that $M$ is orientable. Then there is a
map $f\colon M \to P$ of degree $1$. (In greater
detail, $M=(N\setminus  D \cup (P\setminus  D)$ where
$D$ is a $3$-disk, and $f\colon M \to P$ maps
$N\setminus  D$ to the disk $D$ in $P$ and is the
identity on $P\setminus D$.) Then $f^*\colon
H^3(P;R) \to H^3(M:R)$ is an isomorphism  for every
coefficient ring (group) $R$. Now, for the detecting
element $u$ above, $f^*(u) \ne 0$, and, therefore,
$\wgt(f^*(u)) = 3$. Thus, $f^*(u)$ is a detecting
element for $M$. 

Now, if $M$ is not orientable, then let $\ov M \to
M$ be its orientable double cover (which also is a
closed $3$-manifold). If $\pi_1(M)$ has odd torsion,
then so does $\pi_1(\ov M)$. Because $\ov M$ is
orientable, the argument above says that
$\cat(\ov M) = 3$. But because $\ov M$ covers $M$,
we know that $\cat(\ov M) \leq \cat(M)$. Therefore,
$\cat(M) = 3$. If, on the other hand, there is a
prime component of $M$ with non-free fundamental group  
having no odd torsion, then this component has a
detecting element in $3$-dimensional
$\ZZ/2$-cohomology. Therefore, $M$ has a detecting
element in $\ZZ/2$-cohomology as well and $\cat(M) = 3$.

Now, if $\pi_1(M)$ has odd torsion, then this occurs
in individual prime components. So, $M$ may not have
a detecting element only if we can write $M = P \connsum Q$,
where $P$ is non-orientable and $Q$ is a prime manifold
having odd torsion.

\end{proof}

\m For completeness, note that $\cat(\Sigma) = 1$
for every simply connected $3$-manifold (= homotopy sphere)
$\Sigma$, and, therefore, every non-zero element 
$u\in H^3(\Sigma)$ is a detecting element. Therefore,
we now have proved \thmref{thm:3mancat} and augmented
it by showing that most closed $3$-manifolds possess
detecting elements. The significance of this will
be apparent in \secref{sec:apps}.

\begin{remark}\label{rem:berstein}
In fact, if we allow local coefficients, then all
$3$-manifolds with non-trivial and non-free fundamental
groups have detecting elements. More specifically,
by \cite{Ber}, $\cat(X) = n = \dim(X)$ if and only if
a certain element $u \in H^1(X;I(\pi))$ has $u^n
\not = 0$ in $H^n(X;I(\pi) \otimes \ldots \otimes I(\pi))$.
Here, $\pi = \pi_1(X)$ and $I(\pi)$ is the augmentation
ideal in the group ring $\ZZ\pi$. Since $u^n$ is a cup
product (with local coefficients), it is a detecting
element.
\end{remark}

\section{Two Applications}\label{sec:apps}
A prime motivating problem in the study of
Lusternik-Schnirelmann category has been the
the \emph{Ganea conjecture}; $\cat(X \times S^n)
= \cat(X) + 1$. We now know that the conjecture
is not true in general, so it is even more
interesting to understand when it \emph{is}
valid. For $3$-manifolds, we have the following.

\begin{cory}\label{ganea}
For every closed $3$-manifold $M$,
$$\cat(M\times S^n) = \cat(M)+1 .$$

\noindent That is, the Ganea conjecture holds for $M$.
\end{cory}

\begin{proof}
First, suppose that $M$ is detectable. Then the equality
follows from the general result \cite[Corollary 2.3]{R3},
but the argument in this case is easy. Let $u \in
H^*(M;R)$ have $\wgt(u) = \cat(M)$ and let $v \in
H^n(S^n;R)$ be non-trivial, where, by the results above,
we can always take $R = \ZZ$ or $R = \ZZ/d$. Let
$\tilde u = p_M^*(u)$ and $\tilde v = p_{S^n}^*(v)$,
where $p_M\colon M \times S^n \to M$ and $p_{S^n}
\colon M \times S^n \to S^n$ are the respective
projections. Clearly, $\tilde u \not = 0$ and
$\tilde v \not = 0$ since the compositions
$$\xymatrix{
M \ar[r]^-{\rm incl} & M \times S^n
\ar[r]^-{p_M} & M \\
S^n \ar[r]^-{\rm incl} & M \times S^n
\ar[r]^-{p_{S^n}} & S^n}
$$
are the respective identity maps. By
\propref{prop:swgtprops} (2), $\wgt(\tilde u) \geq
\wgt(u) = \cat(M)$ and $\wgt(\tilde v) \geq \wgt(v)
= 1$. Then the K\"unneth theorem says that
$0 \not = \tilde u \cup \tilde v \in
H^*(M \times S^n;R)$ and (using
\propref{prop:swgtprops} (3) and the product
inequality $\cat(X \times Y) \leq \cat(X) + \cat(Y)$))
$$\cat(M) + 1 \geq \cat(M \times S^n) \geq
\wgt(\tilde u \cup \tilde v) \geq \wgt(\tilde u) +
\wgt(\tilde v) \geq \cat(M) + 1.$$

\noindent Hence, $ \cat(M \times S^n) = \cat(M) + 1$.

Now, suppose that $M$ is not detectable. Then, by
\thmref{thm:non-free}, the oriented double cover $\ov M$
of $M$ is detectable, and $\cat(\ov M) = 3$. 
Therefore, in view of what we said above, $\cat(\ov M
\times S^n) = 4$. But $\ov M \times S^n$ covers
$M \times S^n$, and so $\cat(M \times S^n) \ge 4$. On
the other hand, 
$$
\cat(M \times S^n) \le \cat(M) + 1 = 4
$$
for general reasons. Thus, $\cat(M \times S^n)=4$. 
\end{proof}

\begin{cory}\label{degree} 
Let $f\colon M \to N$ be a degree $1$ map of oriented
$3$-manifolds. Then $\cat M \ge \cat f = \cat N$.
\end{cory}

\begin{proof}
Let $u\in H^3(N;A)$ be a detecting element for $N$.
(Recall that orientable $3$-manifolds always have
detecting elements.) Since $\deg(f)=1$, we 
conclude that $f^*(u) \ne 0$. So, $\cat(f) \ge \wgt(u)$
by \propref{prop:swgtprops} (2). Thus
$$
\cat(M) \ge \cat(f) \ge \wgt(u) = \cat(N).
$$
Of course, $\cat(f) =\cat(N)$ holds since $\cat(f) \leq
\cat(N)$ for general reasons.
\end{proof}

\begin{cory} Let $f\colon M \to N$ be a degree $1$ map
of oriented $3$-manifolds. If $\pi_1(M)$ is free,
then $\pi_1(N)$ is. 
\end{cory}

\begin{proof}
By \corref{degree}, $\cat(N) \le 2$, and so $\pi_1(N)$
is free by \thmref{thm:non-free}
\end{proof}

\

\end{document}